\documentclass[12pt]{amsart}

\newtheorem{theorem}{Theorem}[section]
\newtheorem{lemma}[theorem]{Lemma}
\newtheorem{corollary}[theorem]{Corollary}
\theoremstyle{definition}   
\newtheorem{definition}{Definition}

\theoremstyle{remark}

\numberwithin{equation}{section}

\usepackage{times}
\usepackage{enumerate}
\usepackage{tfrupee}
\usepackage{tikz}
\usetikzlibrary{chains,fit}
\usetikzlibrary{shapes,snakes}
\usetikzlibrary{graphs}
\usepackage{graphicx,adjustbox}
\usepackage{hyperref}
\usepackage{amsmath,amssymb}

\newcommand{\LT}{\mbox{\rm Lt}}
\newcommand{\supp}{\mbox{\rm supp}}

\usepackage{amscd}
\usepackage{graphicx}
\usepackage[all]{xy}

\title[Transversal Intersection and Sum of Polynomial Ideals]
{Transversal Intersection and Sum of Polynomial Ideals}
\author{
Joydip Saha
\and
Indranath Sengupta
\and
Gaurab Tripathi
}
\date{}

\address{\small \rm  Discipline of Mathematics, IIT Gandhinagar, Palaj, Gandhinagar, 
Gujarat 382355, INDIA.} 
\email{saha.joydip56@gmail.com}
\thanks{The first author is a post-doctoral research fellow under the 
research project EMR/2015/000776 sponsored by the SERB, Government of India.}

\address{\small \rm  Discipline of Mathematics, IIT Gandhinagar, Palaj, Gandhinagar, 
Gujarat 382355, INDIA.}
\email{indranathsg@iitgn.ac.in}
\thanks{The second author is the corresponding author, who is supported by the the 
research project EMR/2015/000776 sponsored by the SERB, Government of India.}

\address{\small \rm Department of Mathematics, Jadavpur University, Kolkata,
WB 700 032, India.} 
\email{gelatinx@gmail.com}
\thanks{The third author thanks CSIR for the Senior Research Fellowship.}

\date{}



\begin{document}
\begin{abstract}
In this paper we derive some conditions for transversal intersection of 
polynomial ideals. We exhibit some examples. Finally, as an application 
of the results proved, we compute the Betti numbers for ideals of the form 
$I_{1}(XY) + J$, where $X$ and $Y$ are matrices and $J$ is the ideal generated 
by the $2\times 2$ minors of the matrix consisting of any two rows of $X$. 
\end{abstract}

\maketitle

\section{Introduction}
Ideals $I$ and $J$ are said to intersect transversally if $I\cap J = IJ$. We 
have observed the interesting fact that for transversally intersecting ideals 
$I$ and $J$ in the polynomial ring $R$, the tensor product of minimal free 
resolutions of $R/I$ and $R/J$ is a minimal free resolution of $R/(I+J)$; see 
\ref{tensorprod}. As a part of part of a bigger study of understanding 
the syzygies of ideals of the form $I+J$, where $I$ and $J$ are both determinantal, we 
were motivated to look for criterion for transversal intersection of polynomial ideals; see 
\ref{transideals}. We have come across some natural classes of ideals in the 
polynomial ring which intersect transversally with the rational normal curves and with the 
determinantal ideals of the form $I_{1}(XY)$; see \ref{transversal_IJ}, \ref{transversal_XY}. 
Let us briefly introduce ideals of the form $I_{1}(XY)$ and their sum with other determinantal 
ideals, which are extremely relevant in the field of algebra and geometry and therefore 
forms the central theme of our study. 
\medskip

Let $K$ be a field and $\{x_{ij}; \, 1\leq i,j \leq n\}$, 
$\{y_{j}; \, 1\leq j \leq n\}$ be indeterminates over $K$; $n\geq 2$. Let 
$R:= K[x_{ij}, y_{j}]$ denote the polynomial algebra over $K$. 
Let $X$ denote an $n\times n$ matrix such that its entries are the variables 
$x_{ij}$ and it is either generic or symmetric generic. 
Let $Y=(y_{j})_{n\times 1}$ be the $n\times 1$ column matrix. 
Let $I_{1}(XY)$ denote the ideal generated by the polynomials $g_{j}$, 
which are the $1\times 1$ minors or entries of the $n\times n$ matrix $XY$. 
The primality, primary decomposition 
and Betti numbers of ideals of the form $I_{1}(XY)$ have been studied in \cite{sst} 
and \cite{sstprime}, with the help of Gr\"{o}bner bases for $I_{1}(XY)$.
\medskip
 
Ideals of the form $I_{1}(XY)+J$ are particularly interesting when $J$ is also 
determinantal. They occur in several geometric considerations like linkage and 
generic residual intersection of polynomial ideals, especially in the context 
of syzygies. Bruns-Kustin-Miller \cite{bkm} resolved the ideal 
$I_{1}(XY)+ I_{\min (m,n)}(X)$, where $X$ is a generic $m \times n$ matrix 
and $Y$ is a generic $n\times 1$ matrix. Johnson-McLoud \cite{johnson} proved 
certain properties for the ideals of the form $I_{1}(XY)+I_{2}(X)$, where $X$ 
is a generic symmetric matrix and $Y$ is either generic or generic alternating. 
These ideals 
We have considered the ideal $I_{1}(XY)+I_{2}(\widetilde{X}_{ij})$ (see section 4), 
where $\widetilde{X}_{ij}$ is the matrix consisting of the $i$-th and the $j$-th 
rows of $X$. In an attempt to prove the Cohen-Macaulay property of the 
ring $R/(I_{1}(XY)+I_{2}(\widetilde{X}_{ij}))$, $1\leq i<j\leq n$, we ended 
up with an explicit construction of the minimal free resolution. In the process 
of doing so, we have encountered several examples of transversal intersection of 
ideals, vide Lemmas \ref{transversal} and \ref{transint}; linear quotients, vide 
Lemma \ref{colon} and use the technique of iterated mapping cone along with 
\ref{tensorprod} as an effective tool. 
\bigskip

\section{Criterion for Transversal Intersection of polynomial ideals}
\begin{definition}
Two ideals $I$ and $J$ in the polynomial ring are said to \textit{intersect transversally} 
if $I\cap J = IJ$.
\end{definition}
\medskip

\begin{definition}Let $T\subset R$ be a set of monomials. We define 
\begin{eqnarray*}
{\rm{supp}}(T)& = & \{(i,j,0)\mid x_{ij}\ {\rm divides}\ m \ \text{for\, some}\ m\in T\}\cup \\
{} & {} & \{(0,0,k)\mid y_{k}\ {\rm divides}\ m\ \text{for\, some}\ m\in T\}.
\end{eqnarray*}
If $T = \{m\}$, then we write $\mbox{\rm{supp}}(m)$ instead of $\mbox{\rm{supp}}(\{m\})$.
\end{definition}
\medskip

\begin{lemma} \label{trans}
Let $>$ be a monomial ordering on $R$. Let $I$ and $J$ be ideals in $R$ and let 
$m(I)$ and $m(J)$ denote unique minimal generating sets for their leading ideals 
$\LT(I)$ and $\LT(J)$ respectively. Then, $I\cap J = IJ$ if 
${\rm supp}(m(I)) \cap {\rm supp}(m(J)) = \emptyset$. In other words, 
the ideals $I$ and $J$ intersect transversally if the set of variables 
occurring in the set $m(I)$ is disjointed from the the set of variables 
occurring in the set $m(J)$.
\end{lemma}

\proof Let $f\in (I\cap J)\setminus IJ$. Let $r$ denote 
the remainder term after division of $f$ by a Gr\"{o}bner 
basis of $IJ$ with respect to the monomial order $<$. 
Now $r\in I\cap J$ implies that $r\in I$ and therefore 
${\LT}(r)\in {\LT}(I)$. Hence, there exists monomial 
$m_{i}\in m(I)$ such that $m_{i}\mid {\LT}(r)$. Similarly, 
there exists monomial $m_{j}\in m(J)$ such that 
$m_{j}\mid {\LT}(r)$. Given that $m_{i}$ and $m_{j}$ 
are of disjoint support, we have $m_{i} m_{j}\mid {\LT}(r)$ 
and this proves that ${\LT}(r)\in {\LT}(IJ)$, which is a contradiction. \qed
\medskip

The notion of transversal intersection of ideals $I$ and $J$ become 
particularly useful while resolving ideals of the form $I+J$. 
We say that $I$ and $J$ intersect \textit{transversally} if $I\cap J= IJ$. 
Suppose that $\mathbb{F}_{\cdot}$ resolves $R/I$ and $\mathbb{G}_{\cdot}$ 
resolves $R/J$ minimally. It is interesting to note that if $I$ and $J$ 
intersect transversally, then the tensor product complex 
$\mathbb{F} \otimes_{R} \mathbb{G}$ resolves $R/(I+J)$ minimally; see 
Lemma \ref{tensorprod}. Therefore, it is useful to know if two ideals 
intersect transversally, especially when one is trying to compute minimal 
free resolutions and Betti numbers for ideals of the form $I+J$, through 
iterated techniques; see \cite{hip}.
\medskip

\begin{lemma}\label{tensorprod}
Let $I$ and $J$ be graded ideals in the standard graded polynomial 
ring $R = k[x_{1},\ldots,x_{n}]$ over a field $k$. Let us 
assume that $I\cap J=I\cdot J$. Suppose that $\mathbb{F}_{\centerdot}$ 
and $\mathbb{G}_{\centerdot}$ are minimal graded free resolutions of $I$ and 
$J$ respectively. Then $\mathbb{F}_{\centerdot}\otimes \mathbb{G}_{\centerdot}$ 
is a minimal graded free resolution for the graded ideal $I+J$.
\end{lemma}

\proof Suppose that $\widehat{R}=k[[x_{1},\ldots,x_{n}]]$. 
We have $I\otimes_{R}\widehat{R}\cong I\widehat{R}$ and 
$(I\cap J)\otimes_{R}\widehat{R}\cong I\widehat{R}\cap J\widehat{R}$, since $\widehat{R}$ is a flat 
$R$ algebra (see Theorem 7.4 in \cite{matsumura}). Hence, $I\widehat{R}\cap J\widehat{R}=(I\widehat{R})(J\widehat{R})$. Let $\mathbb{F}_{I}$ and  $\mathbb{F}_{J}$ denote minimal graded free $R$-resolutions of the ideals $I$ and $J$ respectively. Since $\widehat{R}$ is a flat 
$R$ algebra and entries of each matrix that occurs in  $\mathbb{F}_{I}$ 
are homogeneous, we have $\mathbb{E}_{I}=\mathbb{F}_{I}\otimes_{R}\widehat{R}$ 
is a minimal free resolution of $I\widehat{R}$. Similarly, 
$\mathbb{E}_{J}=\mathbb{F}_{J}\otimes_{R}\widehat{R}$ is a minimal 
free resolution of $J\widehat{R}$. 
\medskip

We first prove that $\mathbb{E}_{I}\otimes_{\widehat{R}} \mathbb{E}_{J}$ is a 
minimal free resolution for $I\widehat{R}+J\widehat{R}$. Consider the short exact sequence $0\longrightarrow I\widehat{R} \longrightarrow \widehat{R}
\longrightarrow \widehat{R}/I\widehat{R}\longrightarrow 0$ and tensor it with 
$\widehat{R}/J\widehat{R}$ over $\widehat{R}$. We get the exact sequence 
$$0\longrightarrow {\rm Tor}_{1}^{\widehat{R}}\left( \widehat{R}/I\widehat{R}, \widehat{R}/J\widehat{R}\right) \longrightarrow I\widehat{R}/(I\widehat{R}\cdot J\widehat{R})\longrightarrow \widehat{R}/J\widehat{R}\longrightarrow \widehat{R}/(I\widehat{R}+J\widehat{R}) \longrightarrow 0.$$ 
The terms on the left are $0$ since $\widehat{R}$ is a flat $\widehat{R}$ module. Moreover, the kernel of the map from $I\widehat{R}/I\widehat{R}\cdot J\widehat{R}\longrightarrow \widehat{R}/J\widehat{R}$ is $\widehat{R}\cap J\widehat{R}/I\widehat{R}\cdot J\widehat{R}$. Therefore 
${\rm Tor}_{1}^{\widehat{R}}\left( \widehat{R}/I\widehat{R}, \widehat{R}/J\widehat{R}\right) =0$ 
if and only if $I\widehat{R}\cap J\widehat{R}=I\widehat{R}\cdot J\widehat{R}$. By Corollary 1 
of Theorem 3 proved in \cite{lichtenbaum}, 
${\rm Tor}_{1}^{\widehat{R}}\left( \widehat{R}/I\widehat{R}, \widehat{R}/J\widehat{R}\right) =0$ implies that ${\rm Tor}_{1}^{\widehat{R}}\left( \widehat{R}/I\widehat{R}, \widehat{R}/J\widehat{R}\right) =0$ for all $i\geq 1$. Therefore, $H_{i}(\mathbb{E}_{I}\otimes_{\widehat{R}} \mathbb{E}_{J})\simeq {\rm Tor}_{1}^{\widehat{R}}\left( \widehat{R}/I\widehat{R}, \widehat{R}/J\widehat{R}\right) = 0$ for all $i\geq 1$ and  $H_{0}(\mathbb{E}_{I}\otimes_{\widehat{R}} \mathbb{E}_{J})\simeq \widehat{R}/(I\widehat{R}+J\widehat{R})$. This proves that 
$\mathbb{E}_{I}\otimes_{\widehat{R}} \mathbb{E}_{J}$ resolves $I\widehat{R}+J\widehat{R}$. The resolution is minimal since both $\mathbb{E}_{I}$ and $\mathbb{E}_{J}$ are minimal.
\medskip

We now show that $\mathbb{F}_{I}\otimes_{R}\mathbb{F}_{J}$ is a minimal free resolution of 
$I\widehat{R}+J\widehat{R}$. Let $H_{i}$ be the $i$-th homology of the complex $\mathbb{F}_{I}\otimes_{R}\mathbb{F}_{J}$, then $H_{i}$ is a graded finitely generated $R$-module. Since $\widehat{R}$ is a flat $R$ algebra and 
$\mathbb{E}_{I}\otimes_{\widehat{R}}\mathbb{E}_{J}$ is a minimal free resolution 
of $I\widehat{R}+J\widehat{R}$, we have $H_{i}\otimes_{R}\widehat{R}=0$. Let 
$\mathfrak{m}=\langle x_{1},\ldots,x_{n}\rangle$ be the maximal relevant ideal in the 
standard graded polynomial ring $R$. Now $(H_{i}\otimes_{R}\widehat{R})\otimes_{R}R/m\cong (H_{i}/\mathfrak{m}H_{i}\otimes_{R}\widehat{R}/\mathfrak{m}\widehat{R})\cong (H_{i}/\mathfrak{m}H_{i}\otimes_{R/\mathfrak{m}R}\widehat{R}/\mathfrak{m}\widehat{R})=0$. Therefore 
$H_{i}/\mathfrak{m}H_{i}=0$ and using graded Nakayama $H_{i}=0$. 
Since all entries of matrices that occur in $\mathbb{F}_{I}\otimes_{R}\mathbb{F}_{J}$ 
are homogeneous we have $\mathbb{F}_{I}\otimes_{R}\mathbb{F}_{J}$ 
is a minimal graded free resolution of $I+J$.     \qed
\medskip

\begin{theorem}[Rees]\label{rees} Let $N$ be an $R$ module, $h_{1},\ldots, h_{k}$ 
be an $N$-regular sequence in $R$ and $J=\langle h_{1},\ldots, h_{k}\rangle $. Let $Y=y_{1},\ldots,y_{k}$ be indeterminates over $R$. If $F(y_{1},\ldots ,y_{k})\in N[Y]$ 
is homogeneous of degree $r$ and $F(h_{1},\ldots ,h_{k})\in J^{r+1}N $ then the coefficients of $F$ are in $JN$.
\end{theorem}

\proof See Theorem 1.1.7 in \cite{bh}.\qed
\medskip

\begin{lemma}\label{nzd}
Let $h_{1},\ldots, h_{k}$ be a regular sequence in $R$. Let $J$ denote the ideal 
$\langle h_{1},\ldots,h_{k-1}\rangle$. Then $h_{k}$ is not a zero divisor in $R/J^{r}$, for every 
$r\geq 1$.
\end{lemma}

\proof We use induction on $r$. For $r=1$, the result follows from the 
fact that $h_{1},\ldots, h_{k}$ is a regular sequence in $R$. Let us 
assume that $h_{k}$ is not a zero divisor in $R/J^{r-1}$. Let 
$h_{k}p\in J^{r}$, hence $h_{k}p\in J^{r-1}$. By the induction hypothesis 
we have $p\in J^{r-1}$. We can write 
$$p=\sum_{\lambda_{1}+\cdots+\lambda_{k-1}=r-1}\beta_{(\lambda_{1},\cdots,\lambda_{k-1})}h_{1}^{\lambda_{1}}\cdots h_{k-1}^{\lambda_{k-1}}.$$ 
Let us consider the homogeneous polynomial $F(y_{1},\ldots ,y_{k-1})$ 
of degree $r-1$ in $R[y_{1},\ldots,y_{k-1}]$, given by
$$F(y_{1},\ldots ,y_{k-1})=\sum_{\lambda_{1}+\cdots+\lambda_{k-1}=r-1}h_{k}\beta_{(\lambda_{1},\cdots,\lambda_{k-1})}y_{1}^{\lambda_{1}}\cdots y_{k-1}^{\lambda_{k-1}}.$$ 
Then, $F(h_{1},\ldots ,h_{k-1})=h_{k}p\in J^{r}$. By Theorem \ref{rees}, $\{h_{k}\beta_{(\lambda_{1},\cdots,\lambda_{k-1})}\mid \lambda_{1}+\cdots+\lambda_{k-1}=r-1\}\subset J$. Given that $h_{1},\ldots, h_{k}$ be a regular sequence in $R$, we have $\{\beta_{(\lambda_{1},\cdots,\lambda_{k-1})}\mid \lambda_{1}+\cdots+\lambda_{k-1}=r-1\}\subset J$. Hence, $p\in J^{r}$.\qed
\medskip

\begin{theorem}\label{transideals}
Let $I$ and $J$ be ideals in $R$, such that 
$J$ is generated by an $R/I$ regular sequence $g_{1},g_{2},\cdots, g_{k}$. 
Then,
\begin{enumerate}[(i)]
\item $I\cap J=IJ$.
\item $I\cap J^{r}=IJ^{r}$, for all positive integers $r$.
\end{enumerate}
\end{theorem}

\proof\textbf{(i)} We use induction on the length $k$ of the $R/I$-regular 
sequence generating the ideal $J$. For $k=1$, let $\alpha\in I\cap J$. 
Then since $\alpha\in J$, we can write $\alpha = g_{1}r_{1}$ for some 
$r_{1}\in R$. Therefore, 
$r_{1}g_{1}\in I$ and $\overline{r_{1}g_{1}}= 0\in R/I$. The element 
$\overline{g_{1}}$ is not a zero divisor in $R/I$, by hypothesis. Therefore, 
$\overline{r_{1}}=\overline{0}\in R/I$ and hence $r_{1}\in I$. This shows that $\alpha = r_{1}g_{1}\in IJ$. 
We assume that the statement is true for $k-1$. Let 
$\alpha\in I\cap J$. Since $\alpha\in J$, we can write 
$\alpha = r_{1}g_{1}+r_{2}g_{2}+\cdots+ r_{k}g_{k}$, for some 
$\alpha_{1}, \ldots , \alpha_{k}$ in $R$. Now $\alpha\in I$, therefore $\overline{r_{k}g_{k}}\in R/ (I+ \langle g_{1},g_{2},\cdots, g_{k-1}\rangle)$. 
The elements $\overline{g_{1}},\overline{g_{2}},\cdots, \overline{g_{k}}$ being 
a regular sequence in $R/I$, we have 
$\overline{r}_{k}=0$, that is, $r_{k}\in I+ \langle g_{1},g_{2},\cdots, g_{k-1}\rangle$.
Let $r_{k}= r_{1}^{'}g_{1}+r_{2}^{'}g_{2}+\cdots + r_{k-1}^{'}g_{k-1} + h$, where 
$h\in I$ and $r^{'}_{1},r^{'}_{2},\cdots, r^{'}_{k-1}\in R$. Therefore, 
$\alpha = (r_{1}+r_{1}^{'}g_{k})g_{1}+(r_{2}+r_{2}^{'}g_{k})g_{2}+\cdots + (r_{k-1}+r_{k-1}^{'}g_{k})g_{k-1}+ hg_{k}$. We have $hg_{k}\in I\cap J$. Hence it is 
enough to show that $\alpha^{'}=\alpha- hg_{k}\in IJ$. Let $J^{'}= \langle g_{1},g_{2},\cdots, g_{k-1}\rangle$. Then $\{\overline{g_{1}},\overline{g_{2}},\cdots, \overline{g_{k-1}}\}$ 
being a part of a regular sequence, is a regular sequence in $R/I$. By the induction hypothesis, we have 
$I\cap J'=IJ'$. Now $\alpha\in I$ implies that $\alpha -ig_{k}\in I$ and 
therefore $\alpha'\in I$. Also, $\alpha'\in \langle g_{1},g_{2},\cdots, g_{k-1}\rangle= J^{'}$. Therefore $\alpha'\in I\cap J'= IJ'\subseteq IJ$.\qed
\medskip

\noindent\textbf{(ii)} We use induction on $r$. For $k=1$ the result trivially holds for all $r\geq1$ by (i). We assume by induction that for $I\cap J^{s}=IJ^{s}$ 
for all $1\leq s<r$. Now we 
prove that the result 
holds good for $s=r$. Let $x_{k}\in I\cap J^{r}$. Every element of $J^{r}$ 
can be written in the form  
$$\sum_{i_{1}+i_{2}+\cdots+ i_{k}= r}\alpha_{i_{1},i_{2},\cdots, i_{k}}g_{1}^{i_{1}}\cdots g_{k}^{i_{k}}.$$ 
Therefore, we can write $x_{k} = g_{k}\gamma_{k}+ \beta_{k}$, where 
$\beta_{k}\in \langle g_{1},\cdots, g_{k-1}\rangle^{r}$ and 
$\gamma_{k}\in J^{r-1}$. We know that $g_{k}$ is not a zero divisor in 
$R/(I+\langle g_{1},\cdots, g_{k-1}\rangle )$. It follows by \ref{nzd} that $g_{k}$ 
is also a non-zero divisor in $R/(I+\langle g_{1},\cdots, g_{k-1}\rangle^{r})$. We know that  
$x_{k}\in I\subseteq I+\langle g_{1},\cdots, g_{k-1}\rangle^{r}$. Therefore, 
$\overline{x_{k}}=\overline{0}$ in $R/I+\langle g_{1},\cdots, g_{k-1}\rangle^{r}$ and 
hence $\overline{g}_{k}\overline{\gamma}_{k}=\overline{0}$. This proves that 
$\overline{\gamma}_{k}=\overline{0}$ in $R/I+\langle g_{1},\cdots, g_{k-1}\rangle^{r}$. 
Let $\gamma_{1}= i_{k}+\alpha_{k}$, where $i_{k}\in I$ and $\alpha_{k}\in \langle g_{1},\cdots, g_{k-1}\rangle^{r}$. Hence, $i_{k}=\gamma_{k}-\alpha_{k}\in I\cap J^{r-1}$ and therefore $i_{k}\in I\cap J^{r-1}=IJ^{r-1}$, by the induction hypothesis on $r$. It 
follows that $g_{k}i\in IJ^{r}$. In order to show that $x_{k} = g_{k}(i_{k}+\alpha_{k})+\beta_{k}\in IJ^{r}$, it is therefore enough to prove that 
$x_{k-1}:= g_{k}\alpha_{k} + \beta_{k}\in IJ^{r}$. We have $x_{k-1}\in I \cap \langle g_{1}, \ldots , g_{k-1}\rangle^{r}$. We continue this process to produce $x_{k-2}$, 
$x_{k-3}, \ldots, x_{1}$, such that $x_{i}\in I \cap \langle g_{1}, \ldots , g_{k-i+1}\rangle^{r}$, for every $i=1, \ldots, k$. In particular, $x_{1}\in I\cap \langle g_{1}\rangle^{r} = I\cdot \langle g_{1}\rangle^{r}$. Then, it follows that $x_{2}\in 
I\cdot \langle g_{1},g_{2}\rangle^{r}$. We can successively go back 
and prove that $x_{k}\in IJ^{r}$.\qed
\medskip

Recently, Professor G. Valla pointed out the resemblance of Lemma \ref{transideals} 
in this paper and Lemma 1.1 in \cite{robvalla}.
\bigskip

\section{Transversal Intersection with the Rational Normal Curve} 
Let $S=K[x_{1},\ldots ,x_{n+1}]$, where $x_{i}$'s are indeterminates over 
the field $K$. Let $N =\left(
\begin{array}{ccccc}
x_{1} & x_{2} & x_{3}\cdots &x_{n}\\
x_{2} & x_{3} & x_{4}\cdots &x_{n+1}\\
\end{array}
\right)
$ and let $I_{2}(N)$ denote the ideal generated by the $2\times 2$ minors of the 
matrix $N$ in $S$. The ideal $I_{2}(N)$ is the defining ideal of the rational normal 
curve in the projective space under the standard parametrization. Our aim 
in this section is to show that some natural classes of ideals $J$ in the polynomial 
ring intersect transversally with the ideal $I_{2}(N)$. This information helps us 
write the minimal free resolution of the sum ideal $I_{2}(N)+J$, since the tensor product complex of the minimal free resolutions of $I_{2}(N)$ and $J$ resolve 
$I+J$ minimally; see Lemma 3.7 in \cite{sst}. The main theorems in this section 
are Theorems \ref{transversal_IJ} and \ref{transversal_XY}. We first prove the following 
Lemmas. 
\medskip

\begin{lemma}\label{nzd}
Let $I$ and $J$ be ideals in $S$, such that $J= \langle x \rangle$, where $x$ 
is a non-zero polynomial in $S$. Then $I\cap J = IJ$ if and only if $\overline{x}$ is not a zero 
divisor in $R/I$.
\end{lemma}

\proof Let $\overline{xg}=\overline{0}$ in $R/I$. Then, $xg\in I$ and 
also $xg\in J$. Therefore, $xg\in I\cap J =IJ$. We can write $xg=xg^{'}$, for 
some $g^{'}\in I$. This shows that $x(g-g')=0$ an hence $g=g'\in I$. The 
converse follows from Lemma \ref{transideals}.\qed
\medskip

\begin{lemma}\label{transversal_ab}
The ideals $I_{2}(N)$ and $\langle x_{1}^{a}+x_{n+1}^{b}\rangle $, where $a, b \in \mathbb{N}$, intersects transversally.
\end{lemma} 
 
\proof The ideal $I_{2}(N)$ is kernel of the homomorphism $\zeta:S \rightarrow k[s,t]$ defined as $\zeta(x_{i})=s^{n-i+1}t^{i-1},\quad 1\leq i\leq n+1$. Therefore, 
the ideal $I_{2}(N)$ is a prime ideal. Again, $x_{1}^{a}+x_{n+1}^{b}\notin I_{2}(N)$, since $\zeta( x_{1}^{a}+x_{n+1}^{b})=s^{na}+t^{nb}\neq 0$. Therefore, $\overline{x_{1}^{a}+x_{n+1}^{b}}$ is not a zero divisor in $R/I_{2}(N)$. Hence, by lemma \ref{nzd}, 
the ideals $I_{2}(N)$ and $\langle x_{1}^{a}+x_{n+1}^{b}\rangle$ intersect 
transversally.\qed
\medskip
  
\begin{lemma}\label{nzd3}
Let $a, b$ be natural numbers. Then,
\begin{enumerate}[(i)]
\item The projective dimension of 
$R/(I_{2}(N)+\langle x_{1}^{a}+x_{n+1}^{b}\rangle)$ is $n$.

\item $R/(I_{2}(N)+\langle x_{1}^{a}+x_{n+1}^{b}\rangle)$ is Cohen-Macaulay.
\end{enumerate}
\end{lemma}
 
\proof \textbf{(i)} The minimal free resolution of $R/I_{2}(N)$ 
is given by the Eagon-Northcott complex of length $n-1$. The ideals  
$I$ and $\langle x_{1}^{a}+x_{n+1}^{b}\rangle $ intersect transversally by 
Lemma \ref{transversal_ab}. Therefore, a minimal free resolution of $R/(I_{2}(N)+\langle x_{1}^{a}+x_{n+1}^{b}\rangle)$ is given 
by the tensor product complex of the minimal free resolutions of $R/I_{2}(N)$ 
and $R/\langle x_{1}^{a}+x_{n+1}^{b}\rangle$ by Lemma \ref{tensorprod}. 
Therefore, it follows that the projective dimension of $R/(I_{2}(N)+\langle x_{1}^{a}+x_{n+1}^{b}\rangle)$ is $n(= n-1+1)$.
\medskip

\noindent\textbf{(ii)} Projective dimension of 
$R/(I_{2}(N)+\langle x_{1}^{a}+x_{n+1}^{b}\rangle)$ is $n$, therefore, by 
the Auslander Buchsbaum theorem ${\rm depth}(R/(I_{2}(N)+\langle x_{1}^{a}+x_{n+1}^{b}\rangle))=1 $. Again $I_{2}(N)$ is a prime ideal of height 
$n-1$ and $\overline{x_{1}^{a}+x_{n+1}^{b}}$ is a non-zero divisor 
of $R/I_{2}(N)$, therefore, by the Krull's principal ideal theorem height 
of $I_{2}(N)+\langle x_{1}^{a}+x_{n+1}^{b}\rangle$ is $n$. Hence, 
$\dim(R/(I_{2}(N)+\langle x_{1}^{a}+x_{n+1}^{b}\rangle))=1$ and 
therefore $R/(I_{2}(N)+\langle x_{1}^{a}+x_{n+1}^{b}\rangle)$ is Cohen-Macaulay.\qed
\medskip 
 
\begin{lemma}\label{nzd4}
Let $a,b,c$ be natural numbers. Let $\mathfrak{p}$ be a prime ideal such that 
$I_{2}(N) +\langle x_{1}^{a}+x_{n+1}^{b}, \, x_{n}^{c}\rangle \subset\mathfrak{p}$. 
Then, $\mathfrak{p}=\langle x_{1},\ldots,x_{n+1}\rangle$ and hence the 
height of the ideal $I_{2}(N) +\langle x_{1}^{a}+x_{n+1}^{b},x_{n}^{c}\rangle$ 
is $n+1$.  
\end{lemma}

\proof We use induction to prove $x_{n-i}\in \mathfrak{p}$, for all $0\leq i\leq n-2$. Since $\mathfrak{p}$ is a prime ideal and $x_{n}^{c}\in I_{2}(N) +\langle x_{1}^{a}+x_{n+1}^{b},x_{n}^{c}\rangle \subset\mathfrak{p}$, therefore $x_{n}\in\mathfrak{p}$. Let us assume that $x_{n-i}\in \mathfrak{p}$ for some $i$. Again   
$x_{n-i}x_{n-i-2}-x_{n-i-1}^{2}\in I_{2}(N) +\langle x_{1}^{a}+x_{n+1}^{b},x_{n}^{c}\rangle \subset\mathfrak{p}$. By the induction hypothesis $x_{n-i}\in \mathfrak{p}$, 
therefore we have  $x_{n-i-1}\in \mathfrak{p}$. Therefore  $x_{n-i}\in \mathfrak{p}$ for all $0\leq i\leq n-2$.
\medskip
 
As $x_{1}x_{n+1}-x_{2}x_{n}\in\mathfrak{p}$ and $x_{n}\in\mathfrak{p}$ we have $x_{1}x_{n+1}\in\mathfrak{p}$ hence $ x_{1}$ or $x_{n+1}\in\mathfrak{p}$. If $x_{1}\in \mathfrak{p}$, then $x_{1}^{a}\in \mathfrak{p}$ which implies that $x_{n+1}^{b}\in\mathfrak{p}$. Hence $ x_{n+1}\in \mathfrak{p}$.  If $x_{n+1}\in \mathfrak{p}$, then similarly we can show that $x_{1}\in \mathfrak{p}$ and hence $\mathfrak{p}=\langle x_{1},\ldots,x_{n+1}\rangle$. \qed
\medskip
 
\begin{lemma}\label{nzd5}
$x_{n}^{l}$ is not a zero divisor in $R/(I_{2}(N)+\langle x_{1}^{a}+x_{n+1}^{b}\rangle)$.
\end{lemma}

\proof Suppose that it is a zero divisor, then it is contained in an associated prime ideal of $R/(I_{2}(N)+\langle x_{1}^{a}+x_{n+1}^{b}\rangle)$. But, 
$R/(I_{2}(N)+\langle x_{1}^{a}+x_{n+1}^{b}\rangle)$ being Cohen-Macaulay, the 
prime ideal has to be minimal. We know that the 
height of the ideal $I_{2}(N)+\langle x_{1}^{a}+x_{n+1}^{b}\rangle$ is $n$. Hence, 
any minimal prime ideal of $R/(I_{2}(N)+\langle x_{1}^{a}+x_{n+1}^{b}\rangle)$ 
has height $0$. In other words, any minimal prime ideal containing $(I_{2}(N)+\langle x_{1}^{a}+x_{n+1}^{b}\rangle)$ has height n. But any prime ideal containing 
both $(I_{2}(N)+\langle x_{1}^{a}+x_{n+1}^{b}\rangle)$ and $x_{n}$ has height 
$n+1$ from the previous lemma. Hence it cannot be minimal.\qed 
\medskip

\begin{theorem}\label{transversal_IJ}
Let $J$ denote the ideal $\langle x_{1}^{a}+x_{n+1}^{b}, \, x_{n}^{c}\rangle$ in the polynomial ring $S$, such that $a,b,c$ are in $\mathbb{N}$. Then, $I_{2}(N)\cap J = IJ$. 
\end{theorem}  
 
\proof We have seen in Lemmas \ref{transversal_ab} and \ref{nzd5} that 
$\overline{x_{1}^{a}+x_{n+1}^{b}}, \, \overline{x_{n}^{c}}$ is a regular sequence 
in $R/I_{2}(N)$. Hence, by Theorem \ref{transideals}, the ideals 
$I_{2}(N)$ and $J$ intersect transversally.\qed
\medskip

\begin{theorem}\label{transversal_XY}
Let $T=k[x_{ij},y_{i}]$, where $i,j\in \{1,\ldots,n\}$  Let 
\begin{center}
$X= \begin{pmatrix}
x_{11}&x_{12}&.&.&.&x_{1n} \\
x_{21}&x_{22}&.&.&.&x_{2n} \\
& &. &.& & \\
& &. &.& & \\
& &. &.& & \\
x_{n1}& x_{n2} &.&.&.&x_{nn} \\
\end{pmatrix}$ 
 and  
$Y= \begin{pmatrix}
 y_{1} \\
 y_{2} \\
 . \\
 . \\
 y_{n}
 \end{pmatrix}$
\end{center}
be generic matrices of indeterminates $x_{rs}$ and $y_{k}$.  Let $f_{r} = \sum_{s=1}^{n}x_{rs}y_{s}$, for $1\leq r\leq n$. Then, $I_{1}(XY)=\mathcal{I} = \langle f_{1}, \ldots , f_{n}\rangle$. 
\begin{center}
$$H= \begin{pmatrix}
  x_{i1}&x_{i2}&.&.& x_{i(i-1)}&x_{i(i+1)}&.&.&x_{i(n-1)}&x_{in} \\
 x_{i2}&x_{i3}&.&.& x_{i(i+1)}&x_{i(i+2)}&.&.&x_{in}&x_{pq} \\
 \end{pmatrix},$$
 $1\leq p,q,i\leq n$ and $p\neq i,p\neq q$
 \end{center}
Let $\mathcal{J}=I_{2}(H)$ be rational normal curve. Then $\mathcal{I}\cap \mathcal{J}=\mathcal{I}\mathcal{J}$.
\end{theorem} 

\proof Let us consider the monomial order 
\begin{itemize}
\item $x_{11}> x_{22}> \cdots >x_{nn}$;
\item $x_{rs}, y_{s} < x_{nn}$ for every $1 \leq r \neq s \leq n$. 
\end{itemize}
Then $\{f_{1}, \ldots , f_{n}\}$ forms a Gr\"{o}bner basis of the ideal $\mathcal{I}$ and $\rm{Lt}(f_{r})=x_{rr}y_{r}$. Therefore,  $\supp(\rm{Lt}(\mathcal{I})) \cap \supp(\rm{Lt}(\mathcal{J}))=\emptyset$. Hence by Theorem \ref{dissup} $\mathcal{I}\cap \mathcal{J}=\mathcal{I}\mathcal{J}$. \qed
\bigskip

\section{Resolution of sums of ideals}

\begin{itemize}
\item If $X$ is generic and $i<j$; \, let 
$\widetilde{X}_{ij}=\begin{pmatrix}
x_{i1} & x_{i2} & \cdots & x_{in}\\
x_{j1} & x_{j2} & \cdots & x_{jn}\\
\end{pmatrix}$.

\item If $X$ is generic symmetric and $i<j$; \, let 
$$\widetilde{X}_{ij}=\begin{pmatrix}
x_{1i} & \cdots & x_{ii} & \cdots & x_{ij} & \cdots & x_{in}\\
x_{1j} & \cdots & x_{ij} & \cdots & x_{jj} & \cdots & x_{jn}\\
\end{pmatrix}.$$

\item Let $\mathcal{G}_{ij}$ denote the set of all $2\times 2$ minors of $\widetilde{X}_{ij}$. 

\item Let $I_{2}(\widetilde{X}_{ij})$ denote the ideal generated $\mathcal{G}_{ij}$. 
\end{itemize}
\medskip

\begin{lemma}\label{gb1}
Suppose that $X$ is either generic or generic symmetric. The set $\mathcal{G}_{ij}$ is a Gr\"{o}bner basis for the ideal $I_{2}(\widetilde{X}_{ij})$, with respect a suitable monomial order.
\end{lemma}

\proof We choose the lexicographic monomial order 
given by the following ordering among the variables: $x_{st}> x_{s^{'}t^{'}}$ 
if $(s^{'},t^{'})>_{\rm lex}(s,t)$ and $y_{n}> y_{n-1}>\cdots>y_{1}> x_{st}$ for all $s, t$. 
We now apply Lemma 4.2 in \cite{sst} for the matrix $X^{t}$ and for $k=2$.\qed
\medskip

Our aim in this paper is to prove the following theorem:
\medskip

\begin{theorem}\label{mainthm}
Let $X = (x_{ij})$ be either the generic or the generic symmetric matrix of 
order $n$. Let $1\leq i<j\leq n$. 
\begin{enumerate}
\item The total Betti numbers for the ideal 
$I_{2}(\widetilde{X}_{ij})+\langle g_{i}, g_{j}\rangle$ are given by $b_{0}=1$, $b_{1}= \binom{n}{2}+2$, $b_{2}=2\binom{n}{3}+n$, $b_{i+1}= i\binom{n}{i+1}+(i-2)\binom{n}{i}$, for \, $2\leq i\leq n-2$ \, and \, $b_{n} = n-2$. 
\medskip

\item Let $1\leq k\leq n-2$. Let $\beta_{k,p}$ denote the $p$-th total Betti number for the ideal $I_{2}(\widetilde{X}_{ij}) + \langle g_{i},g_{j},g_{l_{1}}, \ldots , g_{l_{k}}\rangle$, such that $1\leq l_{1} < \ldots < l_{k}\leq n$ and $l_{t}$ is the smallest in the set $\{1,2,\ldots , n\}\setminus \{i,j,l_{1}, \ldots , l_{t-1}\}$, for every $1\leq t\leq k$. They are given by $\beta_{k,0}=1$, $\beta_{k,p} = \beta_{k-1,p-1}+\beta_{k-1,p}$ for $1\leq p\leq n+k-1$ and $\beta_{k,n+k} = n-2$.
\end{enumerate}
\medskip

In particular, the total Betti numbers for the ideal 
$I_{1}(XY)+I_{2}(\widetilde{X}_{ij})$ are $\beta_{n-2,0}, \beta_{n-2,1}, \ldots, \beta_{n-2,2n-2}$. 
\end{theorem}
\bigskip

\section{Preliminaries and some Homological Lemmas}
We first recall some useful results on determinantal ideals pertaining to our work. 
We refer to \cite{concaetal}, \cite{eisenbud}, \cite{peeva} for detailed discussions 
on these. 
\medskip

\begin{lemma}\label{disjoint}
Let $h_{1},h_{2}\cdots, h_{n}\in R$ be such that with respect to a suitable 
monomial order on $R$, the leading terms of them are mutually coprime. 
Then, $h_{1},h_{2}\cdots, h_{n}$ is a regular sequence in $R$.
\end{lemma}

\proof. See Lemma 2.1 in \cite{sstprime}. \qed
\medskip

\begin{theorem}
Let $K$ be a field and let $x_{ij}:1\leq i\leq m,1\leq j\leq n$ be indeterminates over 
$K$. Let $A=(x_{ij})$ be the $m\times n$ matrix of indeterminates and $I_{m}(A)$ denotes 
the ideal generated by the maximal minors of $A$. The set of maximal minors of $A$ is 
a universal Gr\"{o}bner basis for the ideal $I_{m}(A)$.
\end{theorem}

\proof See \cite{concaetal}.\qed
\bigskip

\noindent\textbf{The Eagon-Northcott Complex.} We present the relevant portion from the 
book \cite{eisenbud} here. Let $F=R^{f}$ and $G=R^{g}$ be free modules of finite rank 
over the polynomial ring $R$. The \textit{Eagon-Northcott complex} of a map 
$\alpha: F\longrightarrow G$ (or that of a matrix $A$ representing $\alpha$) 
is a complex 
\begin{eqnarray*}
{EN}(\alpha): 
0\rightarrow (Sym_{f-g}G)^{*}\otimes \wedge^{f}F\stackrel{d_{f-g+1}}\longrightarrow (Sym_{f-g-1}G)^{*}
\otimes \wedge^{f-1}F\stackrel{d_{f-g}} 
\longrightarrow\\
\cdots\longrightarrow(Sym_{2}G)^{*}\otimes
\wedge^{g+2}F\stackrel{d_3}\longrightarrow G^{*}\otimes 
\wedge^{g+1}F\stackrel{d_{2}}\wedge^{g}F\stackrel{\wedge^{g}\alpha}\longrightarrow \wedge^{g} G.
\end{eqnarray*} 
Here $Sym_{k}G$ is the $k$-th symmetric power of G and $M^{*} = {\rm Hom}_{R}(M,R)$. 
The map $d_{j}$ are defined as follows. First we define a diagonal map 
\begin{eqnarray*}
(Sym_{k}G)^{*} & \rightarrow & G^{*}\otimes (Sym_{k-1}G)^{*}\\
u & \mapsto & \sum_{i}u_{i}^{'}\otimes u_{i}^{''}
\end{eqnarray*}
as the dual of the multiplication map $G\otimes Sym_{k-1}G\longrightarrow Sym_{k}G$ 
in the symmetric algebra of $G$. Next we define an analogous diagonal map 
\begin{eqnarray*}
\wedge^{k}F & \longrightarrow & F\otimes \wedge^{k-1}F\\
v & \mapsto & \sum_{i}v_{i}^{'}\otimes v_{i}^{''}
\end{eqnarray*} 
as the dual of the multiplication in the exterior algebra of $F^{*}$.
\medskip

\begin{theorem}[Eagon-Northcott] The Eagon-Northcott complex is a free resolution of $R/I_{g}(\alpha)$ iff grade$(I_{g}(\alpha))=f-g+1$
where $I_{g}(\alpha)$ denotes the $g\times g$ minors of the matrix $A$ representing 
$\alpha$.
\end{theorem}

\proof See \cite{eisenbud}.\qed
\bigskip

\noindent\textbf{Mapping Cone.}\label{MC} We present the relevant portion from the book 
\cite{peeva} here. Let $R$ be the polynomial ring. Let 
$\phi_{\cdot}: (U_{\cdot},d_{\cdot})\rightarrow (U^{'}_{\cdot},d^{'}_{\cdot})$ 
be a map of complexes of finitely generated $R$-modules. The mapping cone of 
$\phi_{\cdot}$ is the complex $W_{\cdot}$ with differential $\delta_{\cdot}$ 
defined as follows. Let $W_{i}=U_{i-1}\oplus U_{i}^{'}$, with 
$\delta|_{U_{i-1}}= -d+ \phi: U_{i-1}\longrightarrow U_{i-2}\oplus U_{i-1}^{'}$ and 
$\delta|_{U_{i}^{'}}= d^{'}: U_{i}^{'}\rightarrow U_{i-1}^{'}$ for each $i$.
\medskip

\begin{theorem}
Let $M$ be an ideal minimally generated by the polynomials $f_{1},\ldots,f_{r}$. 
Set $M_{i}= \langle f_{1},\ldots,f_{i}\rangle$, for $1\leq i\leq r$. 
Thus, $M= M_{r}$. For each $i\geq 1,$ we have the short exact sequence 
$$0\longrightarrow S/(M_{i}:f_{i+1})\stackrel{f_{i+1}}\longrightarrow S/M_{i}\longrightarrow S/M_{i+1}\longrightarrow 0.$$ 
If resolutions of $S/M_{i}$ and $S/(M_{i}:f_{i+1})$ are known then we can 
construct a resolution of $S/M_{i+1}$ by the mapping cone construction.
\end{theorem}

\proof See Construction 27.3 in \cite{peeva}.\qed
\medskip

\begin{lemma}\label{hom1}
Let 
$$R^{a_{1}}\stackrel{A_{1}}{\longrightarrow} R^{a_{2}}\stackrel{A_{2}}{\longrightarrow} R^{a_{3}}$$ 
be an exact sequence of free modules. Let $Q_{1}$, $Q_{2}$, $Q_{3}$ be invertible matrices of sizes 
$a_{1}$, $a_{2}$, $a_{3}$ respectively. Then, 
$$R^{a_{1}}\stackrel{Q_{2}^{-1}A_{1}Q_{1}}{\longrightarrow} R^{a_{2}}\stackrel{Q_{3}^{-1}A_{2}Q_{2}}{\longrightarrow} R^{a_{3}}$$
is also an exact sequence of free modules.
\end{lemma}

\proof The following diagram is a commutative diagram of free modules and the vertical maps are isomorphisms:
$$
\xymatrix{
R^{a_{1}} \ar[r]^{A_{1}}& R^{a_{2}} \ar[r]^{A_{2}} & R^{a_{3}} \\
R^{a_{1}} \ar[u]^{Q_{1}}\ar[r]^{Q_{2}^{-1}A_{1}Q_{1}} & R^{a_{2}} \ar[u]_{Q_{2}}\ar[r]^{Q_{3}^{-1}A_{2}Q_{2}} & R^{a_{3}}\ar[u]_{Q_{3}}
} 
$$
Therefore, 
$R^{a_{1}}\stackrel{Q_{2}^{-1}A_{1}Q_{1}}{\longrightarrow} R^{a_{2}}\stackrel{Q_{3}^{-1}A_{2}Q_{2}}{\longrightarrow} R^{a_{3}}$ 
is exact since 
$R^{a_{1}}\stackrel{A_{1}}{\longrightarrow} R^{a_{2}}\stackrel{A_{2}}{\longrightarrow} R^{a_{3}}$ 
is exact.\qed
\medskip

\begin{corollary}\label{hom2}Let 
$$R^{a_{1}}\stackrel{C}{\longrightarrow} R^{a_{2}}\stackrel{B}{\longrightarrow} R^{a_{3}}\stackrel{A}{\longrightarrow} R^{a_{4}}$$ 
be an exact sequence of free modules. Let $P_{1}$, $P_{2}$, $P_{3}$ be invertible matrices of sizes 
$a_{1}$, $a_{2}$, $a_{3}$ respectively. Then, 
$$R^{a_{1}}\stackrel{P_{2}^{-1}CP_{1}}{\longrightarrow} R^{a_{2}}\stackrel{BP_{2}}{\longrightarrow} R^{a_{3}}\stackrel{AP_{3}^{-1}}{\longrightarrow} R^{a_{4}}$$
is also an exact sequence of free modules.
\end{corollary}

\proof Consider the sequence $R^{a_{1}}\stackrel{C}{\longrightarrow} R^{a_{2}}\stackrel{B}{\longrightarrow} R^{a_{3}}$. If we take $Q_{1}=P_{1},Q_{2}=P_{2}$ and $Q_{3}=I$ and apply Lemma \ref{hom1}, we get that the sequence 
$R^{a_{1}}\stackrel{P_{2}^{-1}CP_{1}}{\longrightarrow} R^{a_{2}}\stackrel{BP_{2}}{\longrightarrow} R^{a_{3}}$ is exact. 
We further note that the entire sequence\, $R^{a_{1}}\stackrel{P_{2}^{-1}CP_{1}}{\longrightarrow} R^{a_{2}}\stackrel{BP_{2}}{\longrightarrow} R^{a_{3}}\stackrel{A}{\longrightarrow} R^{a_{4}}$ is exact as well, since ${\rm Im}(B)= {\rm Im}(BP_{2})$ and $P_{2}$ is invertible. Let us 
now consider the sequence $R^{a_{2}}\stackrel{BP_{2}}{\longrightarrow} R^{a_{3}}\stackrel{A}{\longrightarrow} R^{a_{4}}$. We take 
$Q_{1}=Q_{3}=I$, $Q_{2}= P_{3}^{-1}$ and apply Lemma 2.3 to arrive at our conclusion.\qed
\medskip

\begin{lemma}\label{hom3}
Let 
$$\cdots \longrightarrow R^{\beta_{n+1}}\stackrel{A_{n+1}}{\longrightarrow} R^{\beta_{n}}\stackrel{A_{n}}\longrightarrow R^{\beta_{n-1}}\stackrel{A_{n-1}}\longrightarrow R^{\beta_{n-2}} \longrightarrow \cdots$$ 
be an exact sequence of free R modules. 
Let $a_{ij}$ denote the $(i,j)$-th entry of $A_{n}$. Suppose that $a_{lm}=\pm 1$ for some l and m, 
$a_{li}=0$ for $i\neq m$ and $a_{jm}=0$ for $j\neq l$. Let $A_{n+1}^{'}$ 
be the matrix obtained by deleting the m-th row from $A_{n+1}$, 
$A_{n-1}^{'}$ the matrix obtained by deleting the l-th column from $A_{n-1}$ and 
$A_{n}^{'}$ the matrix obtained by deleting the l-th row and m-th column 
from $A_{n}$. Then, the sequence 
$$\cdots \longrightarrow R^{\beta_{n+1}}\stackrel{A_{n+1}^{'}}{\longrightarrow} R^{\beta_{n}-1}\stackrel{A_{n}^{'}}
\longrightarrow R^{\beta_{n-1}-1}\stackrel{A_{n-1}^{'}}\longrightarrow R^{\beta_{n-2}} \longrightarrow \cdots$$ 
is exact.
\end{lemma}

\proof The fact that the latter sequence is a complex is self evident. We need to prove its exactness. 
By the previous lemma we may assume that $l=m=1$, for we choose elementary matrices to permute rows and 
columns and these matrices are always invertible. Now, due to exactness of the first complex we have 
$A_{n-1}A_{n}=0.$ This implies that the first column of $A_{n-1}$ is $0$, which implies that 
${\rm Im}(A_{n-1}) = {\rm Im}(A_{n-1}^{'})$. Therefore, the right exactness of $A_{n+1}$ is preserved. 
By a similar argument we can prove that the left exactness of $A_{n+1}^{'}$ is preserved.
\medskip

Let $\left({\bf\underline x}\right)$ denote a tuple with entries from $R$. 
If $\left({\bf\underline x}\right)\in {\rm ker}(A_{n}^{'})$, then $\left(0, {\bf\underline x}\right)\in {\rm ker}(A_{n})$. 
There exists $\left({\bf\underline y}\right)\in R^{\beta_{n+1}}$ such that 
$A_{n-1}\left({\bf\underline y}\right) = \left(0, {\bf\underline x}\right)$. It follows that 
$A_{n-1}^{'}\left({\bf\underline y}\right) = \left({\bf\underline x}\right)$, proving the left exactness of 
$A_{n}^{'}$. By a similar argument we can prove the right exactness of $A_{n}^{'}$.\qed
\medskip

\begin{lemma}\label{hom4}
Let $A$ be $q\times p$ matrix over $R$ with $a_{ij}=\pm 1$, 
for some $i$ and $j$. Let $C$ be a $p\times s$ matrix and 
$B$ a $r\times q$ matrix over $R$. There exist an invertible 
$q\times q$ matrix $X$ and an invertible $p\times p$ matrix $Y$, 
such that 
\begin{enumerate}
\item[(i)] $(XAY)_{kj}=\delta_{ki}$ and $(XAY)_{ik}=\delta_{jk}$, 
that is
$$XAY = 
\begin{pmatrix}
\cdots & 0 & \cdots \\
\cdots & 0 & \cdots\\
\vdots & \vdots & \vdots \\
0\cdots 0 & 1 & 0\cdots 0\\
\vdots & \vdots & \vdots \\
\cdots & 0 & \cdots \\
\cdots & 0 & \cdots
\end{pmatrix}; \quad 1 \quad {\rm  at\,\, the} \quad (i,j){\rm -th\,\, spot}.$$ 

\item[(ii)] $(Y^{-1}C)_{kl}= C_{kl}$ for $k\neq j$ and 
$(Y^{-1}C)_{jl}= C_{jl}+ \sum_{t\neq i} (a_{it})C_{tl}$

\item[(iii)] $(BX^{-1})_{kl}= B_{kl}$ for $l\neq i$ and 
$(BX^{-1})_{ki}= B_{ki}+ \sum_{t\neq i} (a_{tj})B_{kt}$.
\end{enumerate}
\end{lemma}

\proof\textbf{(i)} We prove for $a_{ij}=1$. The other case is similar. We take 
$Y= \Pi_{k\neq j}E_{jk}(-a_{ik})$ and $X= \Pi_{k\neq i}E_{ki}(-a_{kj})$, 
where $E_{kl}(\alpha)$ denotes the matrix $E$ with $E_{kl}=\alpha$, $E_{tt}=1$ and 
$E_{ut}=0$ for $u\neq t$ and $(u,t)\neq(k,l)$.
\medskip

\noindent\textbf{(ii)} and \textbf{(iii)} are easy to verify.\qed
\medskip

\begin{lemma}\label{hom5}
Let $A$ be $q\times p$ matrix, $C$ be a $p\times s$ matrix and 
$B$ a $r\times q$ matrix over $R$. The matrices 
$A$, $B$ and $C$ satisfy property $\mathcal{P}_{ij}$ if they satisfy 
the following conditions:
\begin{itemize}
\item $A_{ij}=1$ , $A_{ik}\in \mathfrak{m}$ for $k\neq j$ and $ A_{kj}\in \mathfrak{m}$ for $k\neq i$;
\item $B_{ki}\in\mathfrak{m}$, for $1\leq k\leq r$;
\item $C_{jl}\in\mathfrak{m}$, for $1\leq l\leq s$.
\end{itemize}
The matrices $XAY$, $BX^{-1}$ and $Y^{-1}C$ satisfy property $\mathcal{P}_{ij}$, if 
$A$, $B$, $C$ satisfy property $\mathcal{P}_{ij}$. 
\end{lemma}

\proof This follows from the above lemma since $a_{ik}$ and $a_{kj}$ belong to 
$\mathfrak{m}$.\qed
\bigskip

\section{Betti numbers of $I_{1}(XY)+I_{2}(\widetilde{X}_{ij})$}
\begin{lemma}\label{height}
Let $X$ be generic or generic symmetric matrix. Let $i<j$.
\begin{enumerate}[(i)]
\item ${\rm ht}(I_{2}(\widetilde{X}_{ij}))=n-1$.

\item The Eagon-Northcott complex minimally resolves the ideal $I_{2}(\widetilde{X_{ij}})$.
\end{enumerate}
\end{lemma}

\proof \textbf{(i)} We show that $f_{1}, \ldots , f_{n-1}$, given by 
$f_{k}=x_{ik}x_{j, k+1}-x_{jk}x_{i, k+1}$, $1\leq k \leq n-1$ form a regular sequence.
\medskip

Let us first assume that $X$ is generic. We take the lexicographic monomial order 
induced by the following ordering among the variables: 
$x_{i1}>x_{i2}>\cdots>x_{in}>x_{j2}>x_{j3}>\cdots > x_{jn} > x_{j1}>x_{kl}$, such 
that $x_{kl}$ are those variables which do not appear 
in $\widetilde{X}_{ij}$ and the variables $y_{p}$ are smaller 
than $x_{j1}$. Then, $\LT(f_{k})= x_{ik}x_{j, k+1}$ and hence 
$\gcd(\LT(f_{k}), \LT(f_{l})) = 1$ for every $k\neq l$. Therefore, 
$f_{1}, \ldots , f_{n-1}$ is a regular sequence by Lemma \ref{disjoint} and 
hence ${\rm ht}(I_{2}(\widetilde{X}))\geq n-1$. On the other hand, ${\rm ht}(I_{2}(\widetilde{X})) \leq n-1$, by Theorem [13.10] in \cite{matsumura}. Hence, 
${\rm ht}(I_{2}(\widetilde{X})) \leq n-1$.
\medskip

If $X$ is generic symmetric, then 
we have to choose the lexicographic monomial order induced by 
$x_{ii} > x_{ij} > x_{1i} > x_{2i} > \cdots > x_{i-1,i} > x_{i, i+1} > \cdots > x_{in}>x_{jj}>\cdots > \widehat{x_{ij}} > \cdots > x_{j-1,j} > x_{j,j+1} > \cdots > x_{jn}$ 
and variables $x_{kl}$ not appearing in $\widetilde{X_{ij}}$ and the 
variables $y_{p}$ are smaller than $x_{jn}$.\qed
\medskip

\noindent\textbf{(ii)} The height of $I_{2}(\widetilde{X_{ij}})$ is $n-1$, which is the maximum. Hence, the Eagon-Northcott complex minimally resolves the ideal $I_{2}(\widetilde{X_{ij}})$. \qed
\medskip

\begin{lemma}\label{transversal}
Let $X$ be generic or generic symmetric. Let $i<j$. Then 
$I_{2}(\widetilde{X}_{ij}) \cap \langle g_{i}\rangle 
= I_{2}(\widetilde{X}_{ij}) \cdot \langle g_{i}\rangle$, 
that is, the ideals $I_{2}(\widetilde{X}_{ij})$ and 
$\langle g_{i}\rangle$ intersect transversally.
\end{lemma}

\proof Let $X$ be generic. 
We choose the lexicographic monomial order given by the following ordering among 
the variables: $x_{st}> x_{s^{'}t^{'}}$ if $(s^{'},t^{'})>(s,t)$ and $y_{n}> y_{n-1}>\cdots > y_{1}> x_{st}$ for all $s, t$. Then, by Lemma \ref{gb1} the set of 
all $2\times 2$ minors forms a Gr\"{o}bner basis for the ideal 
$I_{2}(\widetilde{X}_{ij})$. 
Clearly, the minimal generating set $m(\LT(I_{2}(\widetilde{X}_{ij})))$ doesn't involve the indeterminates $x_{in}$ and $y_{n}$, whereas $\LT(g_{i})= x_{in}y_{n}$. Hence, the 
supports of $m(\LT(I_{2}(\widetilde{X}_{ij})))$ and $m(\LT(g_{i}))$ are disjoint. 
Therefore, by Lemma \ref{trans} we are done.
\medskip

Let $X$ be generic symmetric. Once we choose the correct 
monomial order, the rest of the proof is similar to the 
generic case. Suppose that $(i,j)=(n-1,n)$. We choose the lexicographic monomial order given 
by the following ordering among the variables: 
 \begin{eqnarray*}
y_{1}> y_{n}>y_{n-1}>\cdots > y_{2}& > & x_{n-1,n-1} > x_{n-1,n}> x_{1,n-1}> x_{2,n-1}>\cdots > x_{n-2,n-1} \\
{} & > & x_{nn}\\
& > & x_{1n} > \cdots > x_{n-2,n}\\
& > & x_{st} \quad {\rm for\quad all\quad other}\quad s, t.  \qed
\end{eqnarray*}
Suppose that $(i,j)\neq (n-1,n)$. We choose the lexicographic monomial order given 
by the following ordering among the variables: 
\begin{eqnarray*}
y_{n}> y_{n-1}>\cdots > y_{1}& > & x_{ii} >x_{ij}> x_{1i}> x_{2i}>\cdots > x_{i-1,i}>x_{i,i+1}>\cdots > x_{in} \\
{} & > & x_{jj}>\cdots > x_{j-1,j}>x_{j,j+1}>\cdots > x_{jn}\\
& > & x_{st} \quad {\rm for\quad all\quad other}\quad s, t.  \qed
\end{eqnarray*}  
\medskip

\begin{lemma}\label{colon} Let $X$ be generic and $i<j$. Then, 
$( I_{2}(\widetilde{X}_{ij})+\langle g_{i}\rangle: \, g_{j}) = \langle x_{i1}, \ldots, x_{in}\rangle$. If $X$ is generic symmetric and $i<j$, then 
$( I_{2}(\widetilde{X}_{ij})+\langle g_{i}\rangle: \, g_{j}) = \langle x_{1i}, \ldots, x_{i-1,i}, \, x_{ii}, \ldots , x_{in}\rangle$.
\end{lemma}

\proof Let $X$ be generic. We have $x_{it}g_{j}=x_{jt}g_{i}+ \sum_{k=1}^{n}(x_{it}x_{jk}-x_{ik}x_{jt})y_{k}$. Hence, 
$\langle x_{i1}, \cdots, x_{in}\rangle \subseteq \langle I_{2}(\widetilde{X_{ij}}) + \langle g_{i}\rangle: g_{j}\rangle $. Moreover, 
$I_{2}(\widetilde{X}_{ij}) + \langle g_{i}\rangle \subseteq \langle x_{i1},\cdots, 
x_{in}\rangle$ and $g_{j}\notin\langle x_{i1},\cdots ,x_{in}\rangle$. The ideal 
$\langle x_{i1}, \cdots, x_{in}\rangle$ being a prime ideal, 
it follows that $\langle x_{i1},\cdots, x_{in}\rangle \supseteq 
(I_{2}(\widetilde{X}_{ij}) + \langle g_{i}\rangle : g_{j})$. The proof for the generic symmetric case is similar.\qed
\bigskip

\section{Minimal free resolution of 
$I_{2}(\widetilde{X}_{ij})+\langle g_{i},g_{j}\rangle$}
Our aim is to construct a minimal free resolution for the ideal 
$I_{2}(\widetilde{X}_{ij}) + \langle g_{1}, \ldots , g_{n}\rangle$. We have proved 
that the ideals $I_{2}(\widetilde{X}_{ij})$ and $\langle g_{i}\rangle$ intersect transversally; see \ref{transversal}. The ideal 
$I_{2}(\widetilde{X}_{ij}) + \langle g_{i}\rangle$ can therefore 
be resolved minimally by Theorem \ref{tensorprod}. We have also proved 
that the ideal $I_{2}(\widetilde{X}_{ij}) + \langle g_{i}\rangle$ and the ideal 
$\langle g_{j}\rangle$ have linear quotient; see \ref{colon}. Therefore, the 
ideal $I_{2}(\widetilde{X}_{ij}) + \langle g_{i},g_{j}\rangle$ can be resolved by 
the mapping cone construction. A minimal free resolution can then be extracted 
from this resolution by applying Lemma \ref{hom5}. Next, we will show that the 
ideal $I_{2}(\widetilde{X}_{ij}) + \langle g_{i},g_{j}\rangle$ intersects 
transversally with the ideal $\langle g_{l_{1}}\rangle$, if $l_{1}$ is the 
minimum in the set $\{1,2,\ldots , n\}\setminus \{i,j\}$; see Lemma \ref{transint}. Therefore, the ideal 
$I_{2}(\widetilde{X}_{ij}) + \langle g_{i},g_{j},g_{l_{1}}\rangle$ can be 
resolved minimally by Theorem \ref{tensorprod}. Proceeding in this manner, 
we will be able to show that the ideals 
$I_{2}(\widetilde{X}_{ij}) + \langle g_{i},g_{j},g_{l_{1}}, \ldots , g_{l_{k}}\rangle$ 
and $\langle g_{l_{k+1}}\rangle$ intersect transversally, if 
$1\leq l_{1} < \ldots < l_{k} < \l_{k+1}\leq n$ and $l_{k+1}$ is the smallest in the set 
$\{1,2,\ldots , n\}\setminus \{i,j,l_{1}, \ldots , l_{k}\}$; see Lemma \ref{transint}. This finally gives us a minimal free resolution for the ideal $I_{2}(\widetilde{X}_{ij}) + \langle g_{i},g_{j},g_{l_{1}}, \ldots , g_{l_{n-2}}\rangle$, with 
$1\leq l_{1} < \ldots < l_{n-2}\leq n$ and $l_{t}\notin\{i,j\}$ for every $t$. 
\medskip

Let us assume that $X$ is generic and $i=1$ and $j=2$. The proofs for the general 
$i$ and $j$, with $i<j$ would be similar according to the aforesaid scheme. The 
proofs in the case when $X$ is generic symmetric would be similar as well. Comments 
for general $i<j$ and the symmetric case have been made whenever necessary.
\medskip

\subsection{A minimal free resolution for $I_{2}(\widetilde{X}_{12})+\langle g_{1},g_{2}\rangle$}
The minimal free resolution of $I_{2}(\widetilde{X}_{12})$ is given by the Eagon-Northcott 
complex, which is the following:
$$\mathbb{E}_{\centerdot}: 0 \longrightarrow E_{n-1}\longrightarrow\cdots \longrightarrow E_{k}\stackrel{\delta_{k}}\longrightarrow E_{k-1}\longrightarrow \cdots \longrightarrow E_{0}\longrightarrow R/I_{2}(\widetilde{X})\longrightarrow 0 $$
where $E_{0}\cong R^{1}$, $E_{k}= R^{k \binom{n}{k+1}}$ and for each 
$k=0, \ldots , n-2$, the map $\delta_{k} : E_{k}\rightarrow E_{k-1}$ is defined as 
\begin{eqnarray*}
\delta_{k}\left((e_{i_{1}}\wedge\cdots\wedge e_{i_{k+1}})\otimes v_{2}^{k-1}\right) & = 
& \sum_{s=1}^{k+1}x_{2s}(e_{i_{1}}\wedge\cdots\wedge\hat{e_{s}}\wedge\cdots\wedge e_{i_{k+1}})\otimes v_{2}^{k-2}\\
\delta_{k}\left((e_{i_{1}}\wedge\cdots\wedge e_{i_{k+1}})\otimes v_{1}^{k-1}\right) 
 & = & \sum_{s=1}^{k+1}(-1)^{s+1}x_{1s}(e_{i_{1}}\wedge\cdots\wedge\hat{e_{s}}\wedge\cdots\wedge e_{i_{k+1}})\otimes v_{1}^{k-2}\\
\delta_{k}\left((e_{i_{1}}\wedge\cdots\wedge e_{i_{k+1}}\right)\otimes v_{1}^{j}v_{2}^{k-j-1}) & = &\, \sum_{s=1}^{k+1}(-1)^{s+1}x_{1s}(e_{i_{1}}\wedge\cdots\wedge\hat{e_{s}}\wedge\cdots\wedge e_{i_{k+1}})\otimes v_{1}^{j-1}v_{2}^{k-j-1}\\ 
{} & {} & + \, \sum_{s=1}^{k+1}x_{2s}(e_{i_{1}}\wedge\cdots\wedge\hat{e_{s}}\wedge\cdots\wedge e_{i_{k+1}})\otimes v_{1}^{j}v_{2}^{k-j-2}
\end{eqnarray*}
for every ordered $k+1$ tuple $(i_{1},i_{2},\cdots, i_{k+1})$, with 
$1\leq i_{1} < \cdots < i_{k+1}\leq n$ and for every $j=1,2,\cdots, k-2$.
\medskip

A minimal resolution of $\langle g_{1}\rangle$ is given by
$$\mathbb{G}_{\centerdot}: 0\longrightarrow R\stackrel{g_{1}}\longrightarrow R \longrightarrow R/\langle g_{1}\rangle\longrightarrow 0.$$
The ideals $I_{2}(\widetilde{X}_{12})$ and $\langle g_{1}\rangle$ 
intersect transversally, by Lemma \ref{transversal}. Therefore, by Lemma \ref{tensorprod},  
a minimal free resolution for $I_{2}(\widetilde{X}_{12})+\langle g_{1}\rangle$ 
is given by the tensor product complex 
$$\mathbb{E}_{\centerdot}\otimes \mathbb{G}_{\centerdot}: 0\rightarrow E_{n-1}\rightarrow \cdots\rightarrow
E_{k+1}\oplus E_{k}\stackrel{\psi_{k+1}}\longrightarrow E_{k}\oplus E_{k-1}
\rightarrow\cdots \rightarrow E_{0}\rightarrow 
R/I_{2}(\widetilde{X}_{12})+\langle g_{1}\rangle\rightarrow 0$$
such that $\psi_{k}: E_{k}\oplus E_{k-1}\longrightarrow E_{k-1}\oplus E_{k-2}$ is the 
map defined as 
\begin{eqnarray*}
\psi_{k}\left(e_{i_{1}}\wedge\cdots\wedge e_{i_{k+1}}\right)
\otimes v_{1}^{j}v_{2}^{k-j-1} & = &\delta_{k}\left((e_{i_{1}}\wedge\cdots\wedge e_{i_{k+1}})\otimes v_{1}^{j}v_{2}^{k-j-1}\right)\\
\psi_{k}\left(e_{i_{1}}\wedge\cdots\wedge e_{i_{k}}\right)\otimes v_{1}^{j}v_{2}^{k-j-2} & = & (-1)^{k-1}g_{1}\left((e_{i_{1}}\wedge\cdots\wedge e_{i_{k}})\otimes v_{1}^{j}v_{2}^{k-j-2}\right)\\
{} & {} & + \, \delta_{k-1}\left((e_{i_{1}}\wedge\cdots\wedge e_{i_{k}})\otimes v_{1}^{j}v_{2}^{k-j-2}\right).
\end{eqnarray*}

Now we find a minimal free resolution for 
$I_{2}(\widetilde{X}_{12})+\langle g_{1},g_{2}\rangle$ by mapping cone. 
Let $C_{k}:= (\mathbb{E}_{\cdot}\otimes\mathbb{G}_{\cdot})_{k}$. We 
have proved in Lemma \ref{colon} that 
$\langle I_{2}(\widetilde{X}_{12})+\langle g_{1}\rangle: g_{2}\rangle =\langle x_{11},x_{12},\cdots, x_{1n}\rangle$; which is minimally resolved by the Koszul complex. 
Let us denote the Koszul Complex by $(\mathbb{F}_{\cdot}; \sigma_{k})$, 
where $\sigma_{k}$ is the $k$-th differential. We first construct the 
connecting map $\tau_{\cdot}: \mathbb{F}_{\cdot}\rightarrow \mathbb{E}_{\cdot}\otimes\mathbb{G}_{\cdot}$. Let us write $F_{k}:=R^{\binom {n}{k}}$ and $C_{k}:=R^{k\binom {n}{k+1}}\oplus R^{(k-1)\binom {n}{k}}$. The map $\tau_{k}: F_{k}\rightarrow C_{k}$ 
is defined as:
$$\tau_{k}\left(e_{i_{1}}\wedge\cdots\wedge e_{i_{k}}\right) = \sum_{j}y_{j}(e_{i_{1}}\wedge\cdots\wedge e_{i_{k}}\wedge e_{j})\otimes v_{1}^{k-1} - (e_{i_{1}}\wedge\cdots\wedge e_{i_{k}})\otimes v_{1}^{k-2}.$$
Let us choose the lexicographic ordering among the $k$ tuples $(i_{1}, \ldots , i_{k})$, such that $1\leq i_{1} < \cdots < i_{k}\leq n$ 
in order to write an ordered basis for $R^{\binom {n}{k}}$. We define 
lexicographic ordering among the tuples $(i_{1}, \ldots , i_{k+1}, k-j, j)$, 
for $j=0, \ldots , k$ and $k=1, \ldots , n-1$ to order the basis elements for 
$R^{k\binom {n}{k+1}}$. Moreover, in the free module $C_{k}=R^{k\binom {n}{k+1}}\oplus R^{(k-1)\binom {n}{k}}$, we order the basis elements in such a way that those for 
$R^{k\binom {n}{k+1}}$ appear first. The matrix representation of $\tau_{k}$ with 
respect to the chosen ordered bases is the following:
$$\begin{pmatrix}
\mathbf{A}_{k\binom{n}{k+1}\times \binom{n}{k}} & \vline & \mathbf{0}_{k\binom{n}{k+1}\times \binom{n}{k}}\\[6mm]
\hline
&  \vline  & \\
\mathbf{0}_{(k-1)\binom{n}{k}\times \binom{n}{k}} & \vline & \begin{pmatrix}
- \mathbf{I}_{\binom{n}{k}\times \binom{n}{k}}\\[5mm]
\hline\\
\mathbf{0}_{(k-2)\binom{n}{k}\times \binom{n}{k}}
\end{pmatrix}
\end{pmatrix}.$$

\begin{theorem}
The following diagram commutes for every $k = 1, \ldots , n-1$:  
$$
\xymatrix{
F_{k} \ar[r]^{\tau_{k}} & C_{k}\\
F_{k+1}\ar[r]^{\tau_{k+1}}\ar[u]^{\sigma_{k+1}}& C_{k+1} \ar[u]_{\psi_{k+1}}  
} 
$$
\end{theorem}

\proof It suffices to prove the statement for a basis element 
$(e_{i_{1}}\wedge\cdots\wedge e_{i_{k+1}})$ of $F_{k+1}$. 
Without loss of generality we consider $(e_{1}\wedge \cdots \wedge e_{k+1})$. We first compute $(\tau_{k}\circ \sigma_{k+1})(e_{1}\wedge \cdots \wedge e_{k+1})$. 
\begin{eqnarray*}
(e_{1}\wedge \cdots\wedge e_{k+1})& \stackrel{\sigma_{k+1}}\longmapsto & \sum_{j=1}^{k+1}(-1)^{j+1}
x_{1j}(e_{1}\wedge \cdots\ \wedge \hat{e_{j}}\wedge \cdots \wedge e_{k+1})\\
{} & \stackrel{\tau_{k}}\longmapsto & \sum_{j=1}^{k+1}(-1)^{j+1}x_{1j}[\sum_{s=1}^{n}y_{s}(e_1\wedge e_2\wedge\cdots\wedge\hat{e_{j}}\wedge\cdots\wedge e_{k+1}\wedge e_{s})\otimes v_{1}^{k-1}]\\
{} & {} & -\sum_{j=1}^{k+1}(-1)^{j+1}x_{1j}(e_1\wedge e_2\wedge\cdots\wedge\hat{e_{j}}\wedge\cdots\wedge e_{k+1})\otimes v_{1}^{k-2}.
\end{eqnarray*} 

We now compute $(\psi_{k+1}\circ \tau_{k+1})(e_{1}\wedge \cdots \wedge e_{k+1})$.
\begin{eqnarray*}
(e_{1}\wedge \cdots\wedge e_{k+1}) & \stackrel{\tau_{k+1}}\longmapsto & 
\sum_{s=1}^{n}y_{s}(e_{1}\wedge \cdots\wedge e_{k+1}\wedge e_{s})\otimes v_{1}^{k}\\
{} & {} & - (e_1\wedge \cdots\wedge e_{k+1})\otimes v_{1}^{k-1}\\
{} & \stackrel{\psi_{k+1}}\longmapsto & \sum_{s=1}^{n}[\sum_{j=1,2,\cdots,k+1,s}(-1)^{j+1}y_{s}x_{1j}(e_1\wedge \cdots\wedge\hat{e_{j}}\wedge\cdots\wedge e_{k+1}\wedge e_{s})\otimes v_{1}^{k-1}\\
{} & {} & - (-1)^{k}g_{1}(e_{1}\wedge \cdots\wedge\hat{e_{j}}\wedge\cdots\wedge e_{k+1}\wedge e_{j})\otimes v_{1}^{k-1}\\
{} & {} & - \sum_{j=1}^{k+1}(-1)^{j+1}x_{1j}(e_{1}\wedge \cdots\wedge\hat{e_{j}}\wedge\cdots\wedge e_{k+1}\wedge e_{j})\otimes v_{1}^{k-2}\\
{} & = & \sum_{j=1}^{k+1}\sum_{s=1}^{n}[x_{1j}y_{s}(-1)^{j+1}(e_{1}\wedge \cdots\ \wedge\hat{e_{j}}\wedge\cdots\wedge e_{k+1}\wedge e_{s} )\otimes v_{1}^{k-1}]\\
{} & {} & + \sum_{s=1}^{n}(-1)^{s+1}y_{s}x_{1s}(e_{1}\wedge \cdots\wedge\hat{e_{s}}\wedge\cdots\wedge e_{k+1}\wedge e_{s})\otimes v_{1}^{k-1}\\
{} & {} & - (-1)^{k}g_{1}(e_{1}\wedge \cdots\wedge\hat{e_{j}}\wedge\cdots\wedge e_{k+1}\wedge e_{j})\otimes v_{1}^{k-1}\\
{} & {} & - \sum_{j=1}^{k+1}(-1)^{j+1}x_{1j}(e_{1}\wedge \cdots\wedge\hat{e_{j}}\wedge\cdots\wedge e_{k+1}\wedge e_{j})\otimes v_{1}^{k-2}\\
{} & = & \sum_{j=1}^{k+1}\sum_{s=1}^{n}[x_{1j}y_{s}(-1)^{j+1}(e_{1}\wedge \cdots\wedge\hat{e_{j}}\wedge\cdots\wedge e_{k+1}\wedge e_{s} )\otimes v_{1}^{k-1}]\\
{} & {} & + \sum_{s=1}^{n}(-1)^{s+1}(-1)^{k+1-s}y_{s}x_{1s}(e_{1}\wedge \cdots\wedge\hat{e_{s}}\wedge\cdots\wedge e_{k+1}\wedge e_{s})\otimes v_{1}^{k-1}\\
{} & {} & - (-1)^{k}g_{1}(e_{1}\wedge \cdots\wedge\hat{e_{j}}\wedge\cdots\wedge e_{k+1}\wedge e_{j})\otimes v_{1}^{k-1}\\
{} & {} & - \sum_{j=1}^{k+1}(-1)^{j+1}x_{1j}(e_{1}\wedge \cdots\wedge\hat{e_{j}}\wedge\cdots\wedge e_{k+1}\wedge e_{j})\otimes v_{1}^{k-2}\\
{} & = & \sum_{j=1}^{k+1}\sum_{s=1}^{n}[x_{1j}y_{s}(-1)^{j+1}(e_{1}\wedge \cdots \wedge\hat{e_{j}}\wedge\cdots\wedge e_{k+1}\wedge e_{s} )\otimes v_{1}^{k-1}]\\
{} & {} & + (-1)^{k}g_{1}(e_{1}\wedge \cdots\wedge\hat{e_{j}}\wedge\cdots\wedge e_{k+1}\wedge e_{j})\otimes v_{1}^{k-1}\\
{} & {} & - (-1)^{k}g_{1}(e_{1}\wedge \cdots\wedge\hat{e_{j}}\wedge\cdots\wedge e_{k+1}\wedge e_{j})\otimes v_{1}^{k-1}\\
{} & {} & - \sum_{j=1}^{k+1}(-1)^{j+1}x_{1j}(e_{1}\wedge e_{2}\wedge\cdots\wedge\hat{e_{j}}\wedge\cdots\wedge e_{k+1}\wedge e_{j})\otimes v_{1}^{k-2}\\
{} & = & \sum_{j=1}^{k+1}\sum_{s=1}^{n}[x_{1j}y_{s}(-1)^{j+1}(e_{1}\wedge \cdots \wedge\hat{e_{j}}\wedge\cdots\wedge e_{k+1}\wedge e_{s} )\otimes v_{1}^{k-1}]\\
{} & {} & - \sum_{j=1}^{k+1}(-1)^{j+1}x_{1j}(e_{1}\wedge \cdots\wedge\hat{e_{j}}\wedge\cdots\wedge e_{k+1}\wedge e_{j})\otimes v_{1}^{k-2}.\qed
\end{eqnarray*}
\medskip

\noindent Hence the mapping cone $\mathbb{M}(\mathbb{E}_{\cdot}\otimes\mathbb{G}_{\cdot};\mathbb{F}_{\cdot})$ gives us the resolution for $I_{2}(\widetilde{X}_{12})+\langle g_{1},g_{2}\rangle$ as described in \ref{MC}. However, this resolution is not minimal. We now 
construct a minimal free resolution from $\mathbb{M}(\mathbb{E}_{\cdot}\otimes\mathbb{G}_{\cdot};\mathbb{F}_{\cdot})$.
\medskip

A free resolution for the ideal 
$I_{2}(\widetilde{X}_{12})+\langle g_{1},g_{2}\rangle$ has been constructed in 3.1, 
which is given by 
$$0\longrightarrow D_{n+2}\stackrel{d_{n+2}}\longrightarrow D_{n+1}\cdots\stackrel{d_{k+1}}\longrightarrow D_{k}\stackrel{d_{k}}\longrightarrow D_{k-1}\cdots\longrightarrow D_{1}\longrightarrow D_{0}\longrightarrow 0,$$
such that $D_{k}= F_{k-1}\oplus C_{k} = R^{\binom {n}{k-1}}\oplus (R^{k\binom{n}{k+1}}\oplus R^{(k-1)\binom{n}{k}})$ and $d_{k}=(-\sigma_{k-1}+\tau_{k-1}, \, \psi_{k})$. Let us recall that the map 
$\psi$ is the differential in the free resolution for $I_{2}(\widetilde{X}_{12})+\langle g_{1}\rangle$, the map $\sigma$ is the differential in the Koszul resolution for 
$\langle x_{11}, x_{12}, \ldots , x_{1n}\rangle$ and $\tau$ is 
the connecting homomorphism between the complexes defined in 3.1. Let us order 
bases for $F_{k-1}$ and $C_{k}$ with respect to the lexicographic ordering. 
Finally we order basis for $D_{k}$ in such a way that the basis elements for $F_{k-1}$ 
appear first, followed by the basis elements for $C_{k}$. Therefore, the 
matrix representation for the differential map $d_{k}$ is given by
$$\begin{pmatrix}
-\mathbf{\sigma}_{k-1} & \vline & \mathbf{0}\\[6mm]
\hline
&  \vline  & \\
\mathbf{\tau}_{k-1} & \vline & \mathbf{\psi}_{k}
\end{pmatrix} = 
\begin{pmatrix}
-\mathbf{\sigma}_{k-1} & \vline & \mathbf{0}\\[6mm]
\hline
&  \vline  & \\
\begin{pmatrix}
\mathbf{A} & \vline & \mathbf{0}\\[6mm]
\hline
&  \vline  & \\
\mathbf{0} & \vline & \begin{pmatrix}
- \mathbf{I}\\[5mm]
\hline\\
\mathbf{0}
\end{pmatrix}
\end{pmatrix} & \vline & \mathbf{\psi}_{k}
\end{pmatrix}.$$
The entries in the 
matrices representing $\sigma_{k-1}$ and $\psi_{k}$ belong to the 
maximal ideal $\langle x_{ij},y_{j}\rangle$, since both are differentials of minimal 
free resolutions. The block matrix $A$ has also elements in the maximal ideal 
$\langle x_{ij},y_{j}\rangle$. The only block 
which has elements outside the maximal ideal $\langle x_{ij},y_{j}\rangle$ is in the identity block appearing in $\tau_{k-1}$. Therefore, it is clear from the matrix 
representation of the map $d_{k}$ that we can apply Lemma \ref{hom5} repeatedly 
to get rid of non-minimality. Hence, we get a minimal free resolution 
and the total Betti numbers for the ideal 
$I_{2}(\widetilde{X}_{12})+\langle g_{1},g_{2}\rangle$ are 
\begin{eqnarray*}
b_{0} & = & 1,\\[2mm]
b_{1} & = & \binom{n}{2}+2,\\[2mm]
b_{2} & = & 2\cdot \binom{n}{3} + n,\\[2mm]
b_{k+1} & = & k\binom{n}{k+1}+ (k-1)\binom{n}{k}+ \binom{n}{k-1}-\binom{n}{k-1}-\binom{n}{k}\\[1mm]
& = & k\binom{n}{k+1}+ (k-2)\binom{n}{k},\quad {\rm for} \quad 2\leq k\leq n-1,\\[2mm]
b_{n} & = & n-2.
\end{eqnarray*}
\medskip

\subsection{A minimal free resolution for $I_{2}(\widetilde{X}_{ij})+ \langle g_{1},\ldots , g_{n}\rangle$}

\begin{lemma}\label{gb2}
Let $\mathcal{G}_{k}= \mathcal{G}_{12}\cup \{g_{1}, g_{2}, \ldots , g_{k}\}$, $1\leq k\leq n$, where $\mathcal{G}_{12}$ is the set of all $2\times 2$ minors of $\widetilde{X}_{12}$ 
defined in the list of notations in section 2. The set $\mathcal{G}_{k}$ is a Gr\"{o}bner 
basis for the ideal $I_{2}(\widetilde{X}_{12})+\langle g_{1},\ldots , g_{k} \rangle$ with 
respect to a suitable monomial order. 
\end{lemma}

\proof We take the lexicographic monomial ordering in $R$ induced 
by the following ordering among the indeterminates:
\begin{eqnarray*}
x_{nn}> \cdots > x_{tt} > \cdots > x_{33} & > & y_{1}>\cdots > y_{n}\\
{} & > & x_{11} > \cdots > x_{1n} > x_{21} > \cdots > x_{2n}\\
{} &  > & x_{st} \quad {\rm for\,\, other}\quad s,t.
\end{eqnarray*} 
Then, we observe that for every $s\geq 3$, 
$\LT(g_{s})$ is coprime with $\LT(g_{t})$ for every $1\leq t \leq k$; $t\neq s$ and also 
coprime with $\LT(h)$ for every $h\in \mathcal{G}$. Moreover, by Lemma \ref{gb1}, 
$\mathcal{G}_{12}$ is a Gr\"{o}bner basis for $I_{2}(\widetilde{X}_{12})$. Therefore, we only 
have to test the $S$-polynomials $S(g_{1},g_{2})$, $S(g_{1},h)$ and $S(g_{2},h)$, 
for $h\in\mathcal{G}$.
\medskip

We can write $S(g_{1},g_{2})= \sum_{k=1}^{n}[12|1k]y_{k}$ and note that 
$\LT([12|1k])\leq \LT(S(g_{1},g_{2})$ for every $1\leq k\leq n$. Hence, 
$S(g_{1},g_{2})\rightarrow_{\mathcal{G}'} 0$. We note that, if $i\neq 1$ then 
the leading terms of $g_{1}$ and $[12|st]$ are mutually coprime and therefore 
$S(g_{1},[12|st])\rightarrow_{\mathcal{G}_{k}} 0$. Next, the expression 
$S(g_{1},[12|1t])=x_{1t}g_{2}+ \sum_{s\neq t}[12|st]y_{s}$ shows that 
$S(g_{1},[12|1t])\rightarrow_{\mathcal{G}_{k}} 0$. Similarly, if $s\neq 1$ then 
the leading terms of $g_{2}$ and $[12|st]$ are mutually coprime and therefore $S(g_{2},[12|st])\rightarrow_{\mathcal{G}_{k}} 0$. The proof for 
$S(g_{2},[12|1t])$ is similar to that of $S(g_{1},[12|st])$. \qed
\medskip

\noindent\textbf{Remark.} The corresponding result for $i<j$ in general would be the 
following:

\begin{lemma}\label{gb3}
Let $\mathcal{G}_{i,j,k}= \mathcal{G}_{ij}\cup \{g_{i}, g_{j}, g_{l_{1}}, \ldots , g_{l_{k-2}}\}$, $1\leq k\leq n$, $1\leq l_{1}<\cdots < l_{k-2}\leq n$ and 
$l_{t}$ is the smallest in the set $\{1,2,\ldots , n\}\setminus \{i,j,l_{1}, \ldots , l_{t-1}\}$; $\mathcal{G}_{ij}$ denotes the set of all $2\times 2$ minors of $\widetilde{X}_{ij}$ defined in the list of notations in section 2. The set $\mathcal{G}_{i,j,k}$ is a 
Gr\"{o}bner basis for the ideal $I_{2}(\widetilde{X}_{ij})+\langle g_{1},\ldots , g_{k} \rangle$ with 
respect to a suitable monomial order. 
\end{lemma}

\proof While proving this statement with $i<j$ arbitrary, we have to choose the following monomial orders. The rest of the proof 
remains similar. 
\medskip

Suppose that $X$ is generic, we choose the lexicographic monomial ordering in $R$ induced by the following ordering among the indeterminates:
\begin{eqnarray*}
x_{nn}> \cdots > \widehat{x_{jj}}>\cdots > \widehat{x_{ii}}>\cdots >x_{11} 
& > & y_{1}>\cdots > y_{n}\\
{} & > & x_{i1} > \cdots > x_{in}\\
{} & > & x_{j1} > \cdots > x_{jn}\\
{} & > & x_{st} \quad {\rm for\,\, all\,\, other}\quad s,t.
\end{eqnarray*} 

If $X$ is generic symmetric, we choose the lexicographic monomial ordering in $R$ induced by the following ordering among the indeterminates:
\begin{eqnarray*}
x_{nn}> \cdots > \widehat{x_{jj}}>\cdots > \widehat{x_{ii}}>\cdots >x_{11}
> y_{1}>\cdots > y_{n}\\
> x_{ii} > x_{ij} > x_{1i} > x_{2i} > \cdots > x_{i-1,i}>x_{i,i+1}>\cdots > x_{in} \\
> x_{jj} > \cdots > x_{j-1,j}>x_{j,j+1}>\cdots > x_{jn}\\
> x_{st} \quad {\rm for\,\, all\,\, other}\quad s,t.\qed
\end{eqnarray*} 
\medskip

\begin{lemma}\label{transint}
The ideals $I_{2}(\widetilde{X}_{12}) + \langle g_{1}, \ldots , g_{k}\rangle$ and 
$\langle g_{k+1}\rangle$ intersect transversally, for every $2\leq k\leq n-1$.
\end{lemma}

\proof 
Suppose not, then, there exists $h_{k+1}\notin I_{2}(\widetilde{X}_{12})+\langle g_{1}, \ldots , g_{k}\rangle$ such that $h_{k+1}g_{k+1}\in I_{2}(\widetilde{X}_{12})+\langle g_{1}, \ldots , g_{k}\rangle$. Let 
us choose the same monomial order on $R$ as defined in Lemma \ref{gb2}. 
Upon division by elements of $\mathcal{G}_{k}$, we may further assume that $\LT(h)\nmid \LT(h_{k+1})$ for every $h\in\mathcal{G}_{k}$, since 
$\mathcal{G}_{k}$ is a Gr\"{o}bner basis for the ideal $I_{2}(\widetilde{X}_{12})+\langle g_{1},\ldots , g_{k} \rangle$ by Lemma \ref{gb2}. On the other hand 
$h_{k+1}g_{k+1}\in I_{2}(\widetilde{X}_{12})+\langle g_{1}, \ldots , g_{k}\rangle$ and therefore $\LT(h)\mid \LT(h_{k+1})$, for some $h\in\mathcal{G}_{k}$, since $\LT(h)$ and $\LT(g_{k+1})$ are 
mutually coprime, - a contradiction.\qed   
\medskip

\noindent\textbf{Remark.} The corresponding result for $i<j$ in general would be the 
following: The ideals $I_{2}(\widetilde{X}_{ij}) + \langle g_{i},g_{j},g_{l_{1}}, \ldots , g_{l_{k}}\rangle$ and $\langle g_{l_{k+1}}\rangle$ 
intersect transversally, if $1\leq l_{1} < \ldots < l_{k} < \l_{k+1}\leq n$ and 
$l_{k+1}$ is the smallest in the set $\{1,2,\ldots , n\}\setminus \{i,j,l_{1}, \ldots , l_{k}\}$, for every $1\leq k\leq n-3$. The proof is essentially the same as above after 
we use the Lemma \ref{gb3}.
\medskip

\noindent\textbf{Proof of Theorem \ref{mainthm}.} Part (1) of the theorem has been 
proved in 5.1. We now prove part (2) under the assumption $i=1$, $j=2$. Let 
the minimal free resolution of $I_{2}(\widetilde{X}_{12})+\langle g_{1},g_{2},\cdots, g_{k}\rangle$ be $(\mathbb{L}_{\centerdot}, \delta_{\centerdot})$. By Lemma 
\ref{transint} and Lemma \ref{tensorprod}, the minimal free resolution of $I_{2}(\widetilde{X}_{12}) + \langle g_{1},\ldots, g_{k+1}\rangle$ is given by the tensor product of $(\mathbb{L}_{\centerdot}, \delta_{\centerdot})$ and $0\longrightarrow R\stackrel{g_{k+1}}\longrightarrow R\longrightarrow 0$, and that is precisely $(\mathbb{K}_{\centerdot}, \Delta_{\centerdot})$, with $K_{p}= L_{p}\oplus L_{p-1}$ and $\Delta_{p}= (\lambda_{p}, (-1)^{p}g_{k+1}+ \lambda_{p-1})$. Let $\beta_{k,p}$, $0\leq p\leq n+k$, denote the $p$-th total Betti number for the ideal $I_{2}(\widetilde{X}_{12})+\langle g_{1},\ldots,g_{k}\rangle$. Then, the total Betti numbers $\beta_{k+1,p}$, $0\leq p\leq n+k+1$ for the ideal $I_{2}(\widetilde{X}_{12})+\langle g_{1},\ldots,g_{k+1}\rangle$ are given by $\beta_{k+1,0}=1$, $\beta_{k+1,p} = \beta_{k,p-1}+\beta_{k,p}$ for 
$1\leq p\leq n+k$ and $\beta_{k+1,n+k+1} = n-2$. The proof for general $i<j$ 
follows similarly according to the strategy discussed in the beginning of section 5.
\medskip

In particular, the total Betti number $\beta_{n-2,p}$ for the ideal $I_{2}(\widetilde{X}_{ij})+\langle g_{1}, \ldots , g_{n}\rangle$ are given by $\beta_{n-2,0} = 1$, $\beta_{n-2,p} = \beta_{n-3,p-1} + \beta_{n-3,p}$ for $1\leq p\leq 2n-3$ and $\beta_{n-2,2n-2} = n-2$.\qed
\medskip

\noindent\textbf{Example.} We show the Betti numbers at each stage for $n=4$ and $n=5$.
 
$$n=4 :\ \begin{array}{ccccccc}
1&6&8&3&&&\\
1&7&14&11&3&&\\
1&8&12&7&2&&\\
1&9&20&19&9&2&\\
1&10&29&39&28&11&2\\
\end{array}$$

$$n=5 :\ \begin{array}{ccccccccc}
1&10&20&5&4&&&\\
1&11&30&25&9&4&&\\
1&12&25&25&14&3&&\\
1&13&37&50&39&17&3&\\
1&14&50&87&89&56&20&3\\
1&15&64&137&176&145&76&23&3
\end{array}$$
\medskip

\begin{theorem}
The ring $R/I_{1}(XY)+I_{2}(\widetilde{X}_{ij})$, $1\leq i<j\leq n$ is Cohen-Macaulay.
\end{theorem}

\proof We know from \ref{mainthm} that the projective dimension of 
$R/I_{1}(XY)+I_{2}(\widetilde{X}_{ij})$ is $2n-2$. We claim that the 
elements of the set $\mathcal{P}\cup\mathcal{Q}$ forms a regular 
sequence, where $\mathcal{P}=\{x_{ik}x_{j, k+1}-x_{jk}x_{i, k+1}\mid 1\leq k \leq n-1 \}$  
and $\mathcal{Q}=\{g_{t}\mid 1\leq t\leq n,\,\, t\neq j\}$. Suppose 
that $X$ is generic and $j<n$. We consider the matrices 
$$\mathfrak{X}_{ij}=\begin{pmatrix}
x_{1i}&\cdots &\hat{x}_{ji}&\cdots& x_{ni}& x_{ji}\\
x_{1j}&\cdots &\hat{x}_{jj}&\cdots & x_{nj}& x_{jj}\\
\vdots& \vdots& \vdots& \vdots&\vdots&\vdots\\

x_{1n}&\cdots &\hat{x}_{jn}&\cdots & x_{nn}& x_{jn}\\
\end{pmatrix},
\mathfrak{Y}_{ij}=\begin{pmatrix}
y_{1}\\
\vdots\\
\hat{y}_{j}\\
\vdots\\
y_{n}\\
y_{j}
\end{pmatrix}
$$
Then we have, $I_{1}(XY)=I_{1}(\mathfrak{X}_{ij}\mathfrak{Y}_{ij})$. We consider the lexicographic monomial order 
\begin{align*}
& x_{nn}>\cdots >\hat{x}_{jj}>\cdots >\hat{x}_{ii}>\cdots x_{11}>\\
& y_{j}>y_{n}>\cdots >\hat{y}_{j}>\cdots y_{1}>\\
& x_{1i}>\cdots \hat{x}_{ji}>\cdots >x_{ni}> x_{ji}>\\
& x_{1j}>\cdots \hat{x}_{jj}>\cdots >x_{nj}> x_{jj}>\mathrm{other\, indeterminates}.
\end{align*}
Then, $\mathrm{Lt}(g_{i})=x_{ji}y_{j}$ and $\mathrm{Lt}(g_{k})=x_{tt}y_{t}$, 
for $1\leq t\leq n$ and $t\neq i,j$. Therefore, 
$\mathrm{Lt}(x_{ik}x_{j, k+1}-x_{jk}x_{i, k+1})=x_{ik}x_{j, k+1}$ 
for $1\leq k \leq n-1$. The set $\mathcal{P}\cup\mathcal{Q}$ 
forms a regular sequence by \ref{disjoint}, since the leading 
terms of the elements are mutually disjoint. The proof is similar 
in the case $j=n$. Similarly one can prove in the case when $X$ is 
generic symmetric.\qed
\bigskip

\bibliographystyle{amsalpha}

\begin{thebibliography}{A}
\bibitem{bh}
{W. Bruns, J. Herzog, {\em Cohen-Macaulay Rings}, Cambridge Studies in Advanced 
Mathematics 39, Cambridge University Press, UK, 1993.
}

\bibitem{bkm}
{W. Bruns, A.R. Kustin, M. Miller, The Resolution of the Generic Residual Intersection of a Complete Intersection, {\em Journal of Algebra} 128 (1990) 214-239.
}

\bibitem{concaetal}
{A. Conca, Emanuela De Negri, Elisa Gorla, Universal Gr\"{o}bner bases for 
Maximal Minors, {\em International Mathematics Research Notices} 11(2015) 3245-3262.
}

\bibitem{eisenbud}
{D. Eisenbud, {\em Geometry of Syzygies}, Springer-Verlag, NY, 2005.
}

\bibitem{hip}
{P. Gimenez, I. Sengupta and H. Srinivasan,
Minimal graded free resolution for monomial curves defined by arithmetic sequences, 
{\em Journal of Algebra} 388 (2013) 294-310.
}

\bibitem{johnson}
{M.R., Johnson, J. McLoud-Mann, On equations defining Veronese Rings, 
{\em Arch. Math. (Basel)} 86(3)(2006) 205-210.
}

\bibitem{lichtenbaum}
{S.Lichtenbaum, On the vanishing of Tor in regular local rings,
{\em Illinois J.Math.} 10: 220- 226,1966.
}

\bibitem{matsumura}
{H. Matsumura, {\em Commutative Ring Theory}, Cambridge University Press, NY, 1986.
}

\bibitem{miller}
{E. Miller,  B. Sturmfels, {\em Combinatorial Commutative Algebra}, Springer, GTM 227, 2005.
}

\bibitem{peeva}
{I. Peeva, {\em Graded Syzygies}, Springer-Verlag London Limited, 2011.
}

\bibitem{bkm}
{W. Bruns, A.R. Kustin, M. Miller, The Resolution of the Generic Residual Intersection of a Complete Intersection, {\em Journal of Algebra} 128 (1990) 214-239.
}

\bibitem{robvalla}
{L. Robbiano \& G. Valla, On Normal Flatness and Normal Torsion-Freeness, 
{\em Journal of Algebra} 43 (1976) 552 - 560.
}

\bibitem{sst}
{J. Saha, I. Sengupta, G. Tripathi, Ideals of the form $I_{1}(XY)$, 
{\em arXiv:1609.02765 [math.AC]} 2016.
}

\bibitem{sstprime}
{J. Saha, I. Sengupta, G. Tripathi, Primary decomposition and normality of certain determinantal ideals, {\em arXiv:1610.00926 [math.AC]} 2016.
}

\end{thebibliography}

\end{document}